\documentclass[a4paper,10pt]{article}

\usepackage{amsthm}

\usepackage{amsmath}

\usepackage{amsfonts}

\usepackage{amssymb}

\usepackage{amsmath,amsthm,amscd,amssymb}

\usepackage{latexsym}

\usepackage{graphicx}

\usepackage[all]{xy}

\theoremstyle{plain}

\newtheorem{thm}{\textbf{Theorem}}

\newtheorem{lem}{\textbf{Lemma}}

\newtheorem{df}{\textbf{Definition}}

\newtheorem{cor}{\textbf{Corollary}}

\newcommand{\A}{\Bbb{A}}

\newcommand{\R}{\Bbb{R}}

\newcommand{\C}{\Bbb{C}}

\newcommand{\HH}{\Bbb{H}}

\newcommand{\Q}{\Bbb{Q}}

\newcommand{\F}{\Bbb{F}}

\newcommand{\G}{\Bbb{G}}

\newcommand{\Z}{\Bbb{Z}}

\newcommand{\T}{\Bbb{T}}

\newcommand{\p}{\frak{p}}

\newcommand{\GSp}{\text{GSp}}

\newcommand{\Sp}{\text{Sp}}

\newcommand{\GO}{\text{GO}}

\newcommand{\GSO}{\text{GSO}}

\newcommand{\Hom}{\text{Hom}}

\newcommand{\tr}{\text{tr}}

\newcommand{\Gal}{\text{Gal}}

\newcommand{\GL}{\text{GL}}

\newcommand{\y}{\hspace{6pt}}

\title{{\bf{Potential level-lowering for $\GSp(4)$}}}

\author{Claus M. Sorensen}

\begin{document}

\date{}

\maketitle

\begin{abstract}
In this article, we explore a beautiful idea of Skinner and Wiles [SW] in the context of $\GSp(4)$ over a totally real field. The main result provides congruences between automorphic forms which are Iwahori-spherical at a certain place $w$, and forms with a tamely ramified prinicpal series at $w$. Thus, after base change to a finite solvable totally real extension, one can often lower the level at $w$. For the proof, we first establish an analogue of the Jacquet-Langlands correspondence, using the stable trace formula. The congruences are then obtained on inner forms, which are compact at infinity modulo the center, and split at all the finite places. The crucial ingredient allowing us to do so, is an important result of Roche [Roc] on types for principal series representations of split reductive groups.
\footnote{{\it{2000 AMS Mathematics Classification}}: 11F33, 11F41, 11F70, 11F80}
\end{abstract}

\section{Introduction}

Throughout, we fix embeddings $i_{\infty}:\bar{\Q}\rightarrow \C$ and $i_{\ell}: \bar{\Q}\rightarrow \bar{\Q}_{\ell}$ for each prime $\ell$.

\medskip

\noindent In Wiles' proof of Fermat's Last Theorem [Wil], it is shown that certain elliptic curves $E$ over $\Q$ are $modular$: For some/every prime $\ell$, the Galois representation
$$
\text{$\rho_{E,\ell}:\Gal(\bar{\Q}/\Q)\rightarrow \GL(V_{\ell}(E))$, $\y$ $V_{\ell}(E)=\Q_{\ell}\otimes_{\Z_{\ell}}\underset{\longleftarrow}{\lim}E[\ell^n]$,}
$$
on the Tate module is isomorphic to $\rho_{f,\ell}$ for some weight two newform $f$ on $\Gamma_0(N)$. Here
$\rho_{f,\ell}$ is the $\ell$-adic representation attached to $f$, constructed by Eichler and Shimura. It can be realized on the Tate module of a certain quotient of the Jacobian of the modular curve $X_0(N)$. By an ingenious trick, for the proof of modularity, one may assume that $\ell=3$ and that the reduction $\bar{\rho}_{E,3}$ is irreducible.
This is done by the $3\leftrightarrow 5$ switch: Suppose $\bar{\rho}_{E,3}$ is reducible. Then the mod $5$ representation $\bar{\rho}_{E,5}$ must be irreducible since there are no Galois invariant subgroups in $E$ of order $15$. One can then construct another elliptic curve $\mathcal{E}$, with the $same$ mod $5$ representation, such that $\bar{\rho}_{\mathcal{E},3}$ $is$ irreducible. Assuming $\bar{\rho}_{E,3}$ is irreducible, the proof of modularity is divided into two steps:
\begin{itemize}
\item \underline{Residual modularity}: $\bar{\rho}_{E,3}$ is modular, \\
\item \underline{Modularity lifting}: $\bar{\rho}_{E,\ell}$ modular $\Longrightarrow$ $\rho_{E,\ell}$ modular.
\end{itemize}
The first step is essentially a theorem of Langlands and Tunnell. The point is that $\bar{\rho}_{E,3}$ has solvable image, so one can use the theory of base change for $\GL_2$ to get the result. The second step is more involved:
We assume $\bar{\rho}_{E,\ell}$ is irreducible, and denote it by $\bar{\rho}$. For an arbitrary finite set of primes $\Sigma$, we introduce the universal deformation ring $R_{\Sigma}$ classifying $all$ deformations of $\bar{\rho}$, which are well-behaved outside $\Sigma$. In particular, such deformations are unramified at primes outside $\Sigma$ where $\bar{\rho}$ is unramified. In addition we have $\T_{\Sigma}$, the universal modular deformation ring. It is essentially a localization of some Hecke algebra, classifying all $modular$ deformations of $\bar{\rho}$ well-behaved outside $\Sigma$. It is well-defined, at least for large $\Sigma$, exactly because we assume $\bar{\rho}$ has $some$ modular deformation. By universality, there is a natural surjection from $R_{\Sigma}$ onto $\T_{\Sigma}$. The key to step two, is to show that this map is an isomorphism:
$$
R_{\Sigma}\overset{\sim}{\longrightarrow} \T_{\Sigma}.
$$
Wiles gave a $numerical$ criterion for this to be an isomorphism: A certain Selmer group and a certain congruence module should have the same size. Using ideas of Ribet [R1] on level-$raising$, one shows that the sizes of these objects increase in exactly the same way when we add a prime to $\Sigma$. Consequently, it suffices to get an isomorphism in the $minimal$ case where $\Sigma=\varnothing$. However, it is not clear at all that $\T_{\varnothing}$ is well-defined: Does $\bar{\rho}$ have a $minimal$ modular deformation? What is needed, is level-$lowering$ congruences. For example, suppose we have
$$
\text{$f \in S_2(\Gamma_0(Np))$, $\y$ $\bar{\rho}_{f,\ell}$ unramified at $p$,}
$$
where $p \nmid N\ell$. Then there is a form $\widetilde{f}\equiv f$ mod $\ell$ of level $N$. The crucial case of this result was established by Ribet by a complicated geometric analysis of certain Shimura curves [R2]. Once $\T_{\varnothing}$ is known to exist, the minimal case of the isomorphism is finally settled by constructing so-called Taylor-Wiles systems.

\medskip

\noindent Several obstacles occur, when one tries to extend the above to higher-dimensional Galois representations.
For one thing, it is unknown how to generalize the level-raising part. In [CHT], the first installment of the proof of the Sato-Tate conjecture for elliptic curves with non-integral $j$-invariant, the expected results are deduced from a conjectural analogue of Ihara's lemma [Iha]. For $\GSp_4$, we offer $some$ results in this direction in [Sor]. However, inspired by ideas of Kisin, Taylor found a way to bypass Ihara's lemma and complete the proof of Sato-Tate in [T].
Secondly, one does not know in general how to extend Ribet's level-lowering argument. For $U(2,1)$ though, there has been some recent progress by Helm [H]. Instead, in the paper [CHT] for example, one adapts the following elegant idea of Skinner and Wiles: The original paper [SW] was written in the context of Hilbert modular forms, but for this introduction let us for simplicity assume $f$ is the form of level $Np$ above. An argument $much$ simpler than Ribet's, and more amenable to generalization, then shows the existence of a finite solvable totally real extension $F$ over $\Q$, in which all the primes above $\ell$ split, such that
\begin{itemize}
\item $\bar{\rho}_{f,\ell}$ remains irreducible when restricted to $\Gal(\bar{\Q}/F)$, \\
\item there is a form $\widetilde{f}\equiv f_F$ (mod $\ell$) over $F$ of level prime-to-$(p)$.
\end{itemize}
Here $f_F$ denotes the base change of $f$ to $F$. This is of course weaker than Ribet's result quoted above, but it has the strong advantage of being easier to generalize to other situations. Moreover, this weaker version is often enough in order to prove modularity results. Indeed, by the theory of base change for $\GL_2$, it is sufficient to obtain modularity over the solvable extension $F$.

\medskip

\noindent In this paper, we explore the Skinner-Wiles trick for $\GSp_4$ over a totally real field $F$. Our main result is Theorem C below. We first move the initial automorphic representation to an inner form $G$, which is compact at infinity and split at $all$ the finite places. This requires $F$ to have $even$ degree over $\Q$. We accomplish this in Theorem B below, as a consequence of the stable trace formula. Then, for $G$, we obtain the desired congruences in Theorem A. Here, an important ingredient is a result of Roche on types for principal series representations. Using Theorem B again, we then move back to $\GSp_4$. Base change, however, is not yet in such a good shape for $\GSp_4$. Conditional results are available though, due to work of Labesse: Basically, if an automorphic representation has at least two Steinberg components, it admits a weak cyclic base change. We are hopeful that Ngo's recent progress on the fundamental lemma might allow one to weaken the Steinberg condition. Deformations of Galois representations, in the context of $\GSp_4$, has been studied previously by other people. For example, Genestier and Tilouine construct Taylor-Wiles systems in [GT]. We hope that our main result will be of some interest in this area. In the remainder of this introduction, we state precisely our three theorems mentioned just above.

\medskip


\noindent For now, let $F$ be an arbitrary totally real field of degree $d=[F:\Q]$, and let
$$
G_0=\GSp(2n)_{/F}.
$$
Take $G$ to be an arbitrary inner form of $G_0$ over $F$ such that $G_{\infty}$ is compact modulo its center.
To fix ideas, one can think of the unitary similitude group:
$$
G(R)=\{g \in \GL_n(D \otimes_{F}R): {^t\bar{g}}g=\mu(g)I\},
$$
for any $F$-algebra $R$. Here $D$ is a totally definite quaternion $F$-algebra. Its center $Z_G$ is canonically isomorphic to $\G_m$. Let us choose a Hecke character $\omega$, and view it as a central character.
Also, we fix an irreducible representation
$$
\xi:G_{\infty}\rightarrow \GL(V_{\xi}).
$$
Its complexification, an algebraic representation of $G_{\C}^d$, is characterized by its highest weight. In our situation, it corresponds to a tuple $a=(a_v)_{v|\infty}$ where
$$
\text{$a_v=(a_{v,1},\ldots,a_{v,n},\widetilde{a}_v)\in \Z^{n+1}$, $\y$ $a_{v,1}\geq \cdots \geq a_{v,n}\geq 0$.}
$$
Our main result for $G$, is the following analogue of Lemma 4.3.3 in [CHT]:

\medskip

\noindent {\bf{Theorem A.}}
{\it{Let $\ell>2n$ be a prime. Let $\pi$ be an automorphic representation of $G$ with central character $\omega$ and infinity type $\xi$. Suppose $w$ is a finite place of $F$, where the group $G$ splits, such that the following two conditions are satisfied:
\begin{itemize}
\item $\pi_w^{I_w}\neq 0$, where $I_w$ is an Iwahori subgroup in $G_w$,\\
\item $N(\p_w)\equiv 1$ (mod $\ell$).
\end{itemize}
For each finite $v\neq w$ fix a compact open subgroup $K_v$ in $G_v$ such that $\pi_v^{K_v}\neq 0$.
Then there exists an automorphic representation $\widetilde{\pi}$ of $G$, with central character $\omega$ and infinity type $\xi$, satisfying the following three conditions:
\begin{itemize}
\item $\widetilde{\pi}\equiv \pi$ (mod $\ell$),\\
\item $\widetilde{\pi}_w=\widetilde{\chi}_1 \times \cdots \times \widetilde{\chi}_n \rtimes \widetilde{\chi}$ for at most tamely ramified characters $\{\widetilde{\chi}_i,\widetilde{\chi}\}$,\\
\item $\widetilde{\pi}_v^{K_v}\neq 0$ for each finite place $v\neq w$.
\end{itemize}
}}
\medskip

\noindent In this theorem, we use the notation from [Tad] for the principal series representations of $G_0$. Moreover, we use the notation $\widetilde{\pi}\equiv \pi$ (mod $\ell$) to signify that the Hecke eigensystems are congruent modulo $\ell$, via the fixed embeddings $i_{\infty}$ and $i_{\ell}$. More concretely, for almost all $v \nmid \ell$ different from $w$, we have the congruence
$$
\text{$\eta_{\widetilde{\pi}}(\phi)\equiv \eta_{\pi}(\phi)$ (mod $\ell$), $\y$ $\phi \in \mathcal{H}_{\Z}(G_v,K_v)$.}
$$
Here $\eta_{\pi}(\phi)$ gives the action of $\phi$ on the $K_v$-invariant of $\pi_v$, which is a one-dimensional subspace when $K_v$ is hyperspecial. This is true for almost all $v$. One crucial ingredient of the proof of Theorem A, is a result of Roche on types for principal series [Roc]. What we use, is the following special case: Let $\chi_i$ and $\chi$ be $n+1$ characters of $\F_w^*$, inflated to tamely ramified characters of $\mathcal{O}_w^*$. Then, for an irreducible admissible representation $\pi$ of $G_w$, we consider the subspace
of $I_{1}$-invariant vectors transforming according to the characters $\{\chi_i,\chi\}$ under the $I$-action.
Then, this space is $nonzero$ if and only if $\pi$ is a subquotient of
$$
\widetilde{\chi}_1 \times \cdots \times \widetilde{\chi}_n \rtimes \widetilde{\chi}
$$
for $some$ quasi-characters $\widetilde{\chi}_i$ extending the $\chi_i$, and some $\widetilde{\chi}$ extending $\chi$.
Now, the trick is to choose the $\chi_i$ wisely such that the above principal series is irreducible for all possible extensions. Here we use an irreducibility criterion due to Tadi$\acute{\text{c}}$, and the assumptions on $\ell$ and $w$. We can even take the $\chi_i$ and $\chi$ to be $trivial$ modulo $\ell$, and this is used later on in the proof in a crucial way.


\medskip

\noindent To get congruences, as in Theorem A, between automorphic forms on $G_0$, we need to be able to transfer automorphic representations from $G$ to $G_0$ and vice versa. That is, we need an analogue in higher rank of the Jacquet-Langlands correspondence for $\GL_2$ and its inner forms ($n=1$). For this, we assume:
$$
\text{$n=2$, $\y$ $d=[F:\Q]$ is even.}
$$
The first condition seems necessary, in order to apply the trace formula machinery currently available. The second condition ensures the existence (and uniqueness) of a totally definite quaternion algebra $D$ over $F$, split at all finite places. As above, we consider the unitary similitude group: For $F$-algebras $R$,
$$
G(R)=\{g \in \GL_2(D \otimes_{F}R): {^t\bar{g}}g=\mu(g)I\}.
$$
Then $G$ is the inner form of $G_0$, which is compact at infinity mod center, and split at $all$ the finite places.
Throughout the paper, we fix an inner twisting
$$
\psi:G \rightarrow G_0.
$$
That is, an isomorphism over $\bar{\Q}$ such that the Galois actions are intertwined up to $G_0$-conjugation.
Our next theorem provides an analogue of the Jacquet-Langlands correspondence for the groups $G$ and $G_0$. To state it precisely, we need to introduce some notation and terminology. First, we fix a central character $\omega$ and an infinity type $\xi$ as above. Correspondingly, we have an $L$-packet
$$
\text{$\Pi_{\xi}=\otimes_{v|\infty}\Pi_{\xi_v}$, $\y$ $\Pi_{\xi_v}=\{\Pi_{\xi_v}^g,\Pi_{\xi_v}^h\}$}
$$
for $G_{0,\infty}$. Here $\Pi_{\xi_v}^g$ and $\Pi_{\xi_v}^h$ are the generic and holomorphic discrete series representations, respectively, of $\GSp_4(\R)$ having the same infinitesimal character, and the same central character, as $\xi_v$. The theorem we are presenting below, only applies to $stable$ $tempered$ automorphic representations. Here, by $tempered$, we mean that almost all local components are tempered. This implies cuspidality by Langlands's theory of Eisenstein series [La3]: The discrete non-cuspidal spectrum is given by residues of Eisenstein series, and those are non-tempered everywhere, see Proposition 4.5.4 in [Lab]. Our definition of $stable$ is a bit ad hoc: We take it to mean $non$-$endoscopic$ and $non$-CAP, according to

\begin{df}
A cuspidal automorphic representation $\Pi$ of $G_0$, or $G$, is said to be endoscopic if and only if there exists two cuspidal automorphic representations $\rho_1$ and $\rho_2$ of $\GL_2$, with the same central character, such that
for almost all $v$:
$$
L_v(s,\Pi,\text{spin})=L_v(s,\rho_1)L_v(s,\rho_2).
$$
\end{df}

\noindent Here the left-hand side is the Euler factor of the degree four $spin$ $L$-function of $\Pi$. The above definition is from [We2]: See the definition of weak endoscopic lift on page 15, and the pertaining remarks. Let us recall the definition of CAP:

\begin{df}
A cuspidal automorphic representation $\Pi$ of $G_0$, or $G$, is said to be CAP (cuspidal associated to parabolic) relative to the parabolic $P$ in $G_0$ if and only if there exists a cuspidal automorphic representations $\rho$ of a Levi factor $M_P$ such that $\Pi$ is weakly equivalent to the constituents of $\text{Ind}_P^{G_0}(\rho)$.
\end{df}

\noindent For $\GSp_4$, these are completely understood due to the work of Piatetski-Shapiro [PS] and Soudry [Sou]. One can construct them by $\theta$-correspondence. Let us mention that Arthur's conjecture predicts that non-CAP $implies$ tempered, but this is not known over $F$. We can now state our main result on functoriality:

\medskip

\noindent {\bf{Theorem B.}}
{\it{Fix an arbitrary member $\Pi_{\xi}^{+}$ of the $L$-packet $\Pi_{\xi}$. Then there is a natural multiplicity-preserving one-to-one correspondence between the sets:
$$
\{\text{stable tempered automorphic $\pi$ of $G$ with $\omega_{\pi}=\omega$ and $\pi_{\infty}=\xi$}\}\overset{1:1}{\leftrightarrow}
$$
$$
\{\text{stable tempered cuspidal automorphic $\Pi$ of $G_0$ with $\omega_{\Pi}=\omega$ and $\Pi_{\infty}=\Pi_{\xi}^{+}$}\}.
$$
The correspondence takes $\pi \mapsto \Pi_{\xi}^{+}\otimes \pi_f$ and $\Pi \mapsto \xi\otimes \Pi_f$.}}
\medskip

\noindent The proof of this theorem is an exercise based on deep results of Arthur, Hales, Langlands, Shelstad, Waldspurger, and Weissauer. In [Sor] we proved an analogous result over $\Q$. However, in that case, the inner form $G$ is necessarily non-split at some finite place. This makes it impossible to move $from$ $G_0$ to $G$, as there are no known character relations a la Shelstad at this finite place.


\medskip

\noindent Combining the two preceding theorems, we obtain congruences for $\GSp_4$ a la Skinner and Wiles [SW] in the case of $\GL_2$. Thus, we continue to assume that $n=2$ and $d$ is even. The strategy is clear: We start out with a $\Pi$ on $G_0$. By Theorem B, it corresponds to a $\pi$ on $G$ with the same finite part. Theorem A then yields a $\widetilde{\pi}\equiv \pi$ (mod $\ell$). Now, we wish to apply Theorem B to $\widetilde{\pi}$ and find a corresponding representation $\widetilde{\Pi}$ on $G_0$. For that we need to make sure $\widetilde{\pi}$ is stable and tempered. Hence, we assume the original $\Pi$ admits a Galois representation $\rho_{\Pi,\ell}$ and that its mod $\ell$ reduction $\bar{\rho}_{\Pi,\ell}$ is irreducible. Since $\widetilde{\pi}$ is congruent to $\pi$ mod $\ell$, this ensures that
$\widetilde{\pi}$ is stable (endoscopic forms and CAP forms have reducible Galois representations). Conjecturally, it is automatically tempered, as mentioned above. One only expects to attach Galois representations to $algebraic$ automorphic representations of $G_0$.

\medskip

\noindent {\bf{Definition 3.}}
{\it{A cuspidal automorphic representation $\Pi$ of $G_0$ is called algebraic if and only if the following holds:
For each archimedean place $v$ of $F$ consider the $L$-parameter of $\Pi_v$, say $\sigma_{\Pi_v}$ [La1]. This is an admissible homomorphism
$$
\sigma_{\Pi_v}:W_{F_v}\rightarrow {^LG_0}=\GSp_4(\C)\times W_{F_v},
$$
where $W_{F_v}$ is the Weil group. When we restrict $\sigma_{\Pi_v}$ to $\C^*$, the corresponding four-dimensional representation should be a sum of characters of the form}}
$$
\text{$z \mapsto z^p\bar{z}^q$, $\y$ {\it{for}} $\y$ $p,q \in \Z$.}
$$
\noindent Algebraicity gives a restriction on the possible infinity types $\xi$. For example, in the Hilbert modular case $n=1$, the classical weights must have the same parity. By a principle propounded by Langlands, algebraic cusp forms on $\GL_n$ should be associated with pure motives. See [Ram] for a precise statement. Here we will only need the $\ell$-adic Galois representation $\rho_{\Pi,\ell}$. We admit the following:

\medskip

\noindent {\bf{Conjecture.}}
{\it{Let $\Pi$ be a cuspidal automorphic representation of $G_0$ over $F$, with infinity type in some $\Pi_{\xi}$, and assume the twist $\Pi\otimes|\det|^{w/2}$ is algebraic for some integer $w$. Let $S_{\Pi}$ denote the set of ramified places of the representation $\Pi$. Then there exists a four-dimensional continuous semisimple Galois representation
$$
\rho_{\Pi,\ell}:\Gal(\bar{\Q}/F)\rightarrow \GL_4(\bar{\Q}_{\ell}),
$$
which is unramified at all $v \notin S_{\Pi}$ not dividing $\ell$. Moreover, for such places $v$,
$$
L_v(s-w/2,\Pi,\text{spin})=\det(1-\rho_{\Pi,\ell}(\text{Frob}_v)N(\p_v)^{-s})^{-1}.
$$
Here $\text{Frob}_v$ is the geometric Frobenius. If $\Pi$ is not CAP, the representation
$\rho_{\Pi,\ell}$ is pure of weight $w$. That is, the eigenvalues of $\rho_{\Pi,\ell}(\text{Frob}_v)$ all have complex absolute value $N(\p_v)^{w/2}$ for $v$ as above. In other words, each $\Pi_v$ is tempered.
}}

\medskip

\noindent This was proved over $\Q$ by Weissauer [W3]. Therefore, if $\Pi$ is a base change from $\Q$, the conjecture holds. In fact, cyclic base change is $known$ for $G_0$ in many cases: Suppose $F$ is a cyclic extension of $\Q$, and let $\Pi$ be a cuspidal automorphic representation of $G_0$ over $\Q$, with at least $two$ Steinberg components. Then $\Pi$ admits a weak cuspidal base change $\Pi_F$ on $G_0$ over $F$. This is due to Labesse in great generality. See Theorem 4.6.2 in [Lab] for a more precise statement.

\medskip

\noindent {\bf{Theorem C.}}
{\it{Let $\ell>3$ be a prime. Let $\Pi$ be a cuspidal automorphic representation of $G_0$ with central character $\omega$ and infinity type $\Pi_{\xi}^{+}$. We assume $\Pi$ admits a Galois representation $\rho_{\Pi,\ell}$ as above, and that its mod $\ell$ reduction $\bar{\rho}_{\Pi,\ell}$ is irreducible. Suppose $w$ is a finite place of $F$, satisfying the two conditions:
\begin{itemize}
\item $\Pi_w^{I_w}\neq 0$, where $I_w$ is an Iwahori subgroup in $G_{0,w}$,\\
\item $N(\p_w)\equiv 1$ (mod $\ell$).
\end{itemize}
For each finite $v\neq w$ fix a compact open subgroup $K_v$ in $G_v$ such that $\Pi_v^{K_v}\neq 0$.
Then there exists a cuspidal automorphic representation $\widetilde{\Pi}$ of $G_0$, with central character $\omega$ and infinity type $\Pi_{\xi}^{+}$, satisfying the following three conditions:
\begin{itemize}
\item $\widetilde{\Pi}\equiv \Pi$ (mod $\ell$),\\
\item $\widetilde{\Pi}_w=\widetilde{\chi}_1 \times \widetilde{\chi}_2 \rtimes \widetilde{\chi}$ for at most tamely ramified characters $\{\widetilde{\chi}_1,\widetilde{\chi}_2,\widetilde{\chi}\}$,\\
\item $\widetilde{\Pi}_v^{K_v}\neq 0$ for each finite place $v\neq w$.
\end{itemize}
}}
\medskip

\noindent What keeps us from getting a precise analogue of the main result in [SW], is that $strong$ cyclic base change is unavailable for $\GSp(4)$. Gan and Takeda have recently shown the local Langlands conjecture in this case [GaT], so locally one can base change $L$-packets. What we need is a global lift compatible with the local lift at $all$ places.
In the following informal discussion, we assume we have such a strong lift: We take $F$ to be $any$ totally real field, and let $\Pi$ be a cuspidal automorphic representation of $G_{0/F}$ with a Galois representation $\rho_{\Pi,\ell}$. Here $\ell>3$ is a prime such that the reduction $\bar{\rho}_{\Pi,\ell}$ is irreducible. First, one can base change to a situation where Theorem C applies. Namely, there is a finite solvable totally real extension $E/F$ of even degree with the following properties:
\begin{itemize}
\item $\bar{\rho}_{\Pi,\ell}|_{\Gal(\bar{\Q}/E)}$ is irreducible,\\
\item $\text{BC}_{E_{\widetilde{w}}/F_w}(\Pi_w)$ is Iwahori-spherical\footnote{To be precise, the base change of the $L$-packet of $\Pi_w$ is an Iwahori-spherical singleton} for all $\widetilde{w}|w$,\\
\item $N(\p_{\widetilde{w}})\equiv 1$ (mod $\ell$) for all $\widetilde{w}|w$.
\end{itemize}
Here $w\nmid \ell$ is a fixed finite place of $F$ such that $\Pi_w$ is ramified. Moreover, to retain some control above $\ell$, we may assume that all the places above $\ell$ split completely in $E$. The existence of $E$ follows from a well-known fact, also used in [SW], that one can always find a totally real cyclic extension with prescribed splitting and ramification at a finite set of places. Now, Theorem C applies to any base change of $\Pi$. One deduces the existence of a cuspidal $\widetilde{\Pi}$ on $G_{0/E}$ with
\begin{itemize}
\item $\omega_{\widetilde{\Pi}}=\omega_{\Pi}\circ N_{E/F}$ and $\widetilde{\Pi}_{\infty}=\Pi_{\infty}^{\otimes [E:F]}$,\\
\item $\bar{\rho}_{\widetilde{\Pi},\ell}\sim \bar{\rho}_{\Pi,\ell}|_{\Gal(\bar{\Q}/E)}$,\\
\item $\widetilde{\Pi}_{\widetilde{w}}$ is a tamely ramified principal series for all $\widetilde{w}|w$.
\end{itemize}
Furhermore, we may arrange that $\widetilde{\Pi}_{\widetilde{v}}$ is unramified at all $\widetilde{v}|v$ such that $\Pi_v$ is unramified. Recall that the places above $\ell$ split in $E$, so we may also arrange for $\widetilde{\Pi}$ to have the $same$ level as $\Pi$ at such places. Finally, after another solvable base change, we may assume 
$\widetilde{\Pi}_{\widetilde{w}}$ is an $unramified$ principal series for all $\widetilde{w}|w$. That is, we can lower the level by passing to a solvable extension. Using this result, an inductive argument shows the existence of a finite solvable extension $E/F$ such that $\bar{\rho}_{\Pi,\ell}|_{\Gal(\bar{\Q}/E)}$ admits an automorphic deformation $\rho_{\widetilde{\Pi},\ell}$ with $minimal$ ramification. However, we stress that our discussion in this paragraph relies heavily on the availability of $strong$ cyclic base change for $\GSp(4)$.

\medskip

\noindent I would like to thank A. Wiles for a question instigating this project. While working on this paper, I have benefited greatly from enlightening conversations and correspondence with J. Getz, D. Ramakrishnan, and C. Skinner.

\section{Principal series representations of $\GSp(2n)$}

In this section we review two results on principal series needed later. In this section only, $F$ denotes a finite extension of $\Q_p$. We use $\nu$ to signify its normalized absolute value $|\cdot|_F$, that is, the modulus character of $F$. We let $\p$ denote the maximal ideal in the valuation ring $\mathcal{O}$, and choose a uniformizer $\varpi$. Also,
$$
G_0=\GSp(2n),
$$
viewed as an algebraic group over $F$. Here the similitude group is defined using
$$
\mathcal{J}=\begin{pmatrix}& & 1 \\ & \ddots & \\ -1 & & \end{pmatrix}.     
$$
The group $G_0$ comes with the similitude character $\mu$. So, for any $F$-algebra $R$,
$$
G_0(R)=\{g \in \GL_{2n}(R): {^tg}\mathcal{J}g=\mu(g)\mathcal{J}\}.
$$
The subset of upper triangular matrices in $G_0$ form a Borel subgroup $B$, and the diagonal matrices form a split maximal torus $T$ of dimension $n+1$. Indeed,
$$
T(R)=\{t=\begin{pmatrix}t_1 & & & & &\\ & \ddots & & & & \\ & & t_n & & & \\ & & & \mu t_n^{-1} & &\\ & & & & \ddots & \\ & & & & & \mu t_1^{-1}\end{pmatrix}\},
$$
for any $R$ as above. Note that $n$ is the rank of the derived group $G_0^{\text{der}}$. To define principal series, consider quasi-characters $\chi_1,\ldots, \chi_n$ and $\chi$ of $F^*$, and look at
$$
\chi_1 \otimes \cdots \otimes \chi_n \otimes \chi: t \mapsto \chi_1(t_1)\cdots \chi_n(t_n)\chi(\mu).
$$
We view it as a character of $B(F)$, and induce it to a representation of $G_0(F)$:
$$
\chi_1 \times \cdots \times \chi_n \rtimes\chi= \text{Ind}_B^{G_0}(\chi_1 \otimes \cdots \otimes \chi_n \otimes \chi).
$$
The induction is normalized by $\delta_B^{\frac{1}{2}}$ in the standard way, where $\delta_B$ is the modulus character of $B(F)$. Here we use the notation introduced by Tadi$\acute{\text{c}}$ in [Tad].

\subsection{An irreducibility criterion of Tadi$\acute{\text{c}}$}

For simplicity, we now assume $p$ is $odd$. Then $F^*$ has exactly three characters of order two, and this simplifies the irreducibility criterion we are about to describe. The reducibility of unitary principal series for symplectic groups was first studied by Winarsky [Win], and then extended by Keys in [Key]. Using their results, Tadi$\acute{\text{c}}$ obtained complete results for similitude groups. We need:

\begin{thm}
The principal series $\chi_1 \times \cdots \times \chi_n \rtimes\chi$ is irreducible if and only if
\begin{itemize}
\item There are at most two distinct $\chi_i$ of order two, \\
\item $\chi_i \neq \nu^{\pm 1}$, for all $i$,\\
\item $\chi_i \neq \nu^{\pm 1}\chi_j^{\pm 1}$, for all $i<j$.
\end{itemize}
\end{thm}

\noindent $Proof$. This is Theorem 7.9 combined with Remark 7.10 in [Tad]. $\square$

\medskip

\noindent In the last condition, all four combinations of signs are allowed. If all the $\chi_i$ are unitary, the last two conditions are clearly satisfied. In particular, for the group $\GSp(4)$, unitary principal series are always irreducible. For arbitrary $n$, the length is at most two - see Corollary 7.8 in [Tad]. Now, let us fix a prime
$$
\text{$\ell>2n$ such that $N(\p)\equiv 1$ (mod $\ell$).}
$$
Here $N(\p)$ is the cardinality of the residue field $\F$. We assume such an $\ell$ exists.

\begin{cor}
There exists $n$ tamely ramified characters $\chi_i$ of $\mathcal{O}^*$ such that
$$
\widetilde{\chi}_1 \times \cdots \times \widetilde{\chi}_n \rtimes\widetilde{\chi}
$$
is irreducible for all quasi-characters $\widetilde{\chi}_i$ with $\widetilde{\chi}_i|_{\mathcal{O}^*}=\chi_i$, and all $\widetilde{\chi}$. Moreover,
$$
\text{$\chi_i \equiv 1$ (mod $\ell$).}
$$
\end{cor}

\noindent $Proof$. We choose a character $\chi$ of $\F^*$ of order $\ell$, and inflate it to a character of the unit group $\mathcal{O}^*$. We then take $\chi_i$ to be $\chi^i$ for each $i$. By assumption on $\ell$,
$$
\{\chi_1,\ldots,\chi_n\} \cap \{\chi_1^{-1},\ldots,\chi_n^{-1}\}=\varnothing.
$$
If $\widetilde{\chi}_i$ are quasi-characters of $F^*$ extending the $\chi_i$, the three irreducibility conditions above are satisfied: None of the tamely ramified $\widetilde{\chi}_i$ have order two,
$$
\widetilde{\chi}_i=\nu^{\pm 1}\widetilde{\chi}_j^{\pm 1}\Longrightarrow \chi_i=\chi_j^{\pm 1}.
$$
This contradicts the fact that the above intersection is empty. To show $\chi_i$ is trivial mod $\ell$, note that its reduction maps into the $\ell$th roots of unity in $\bar{\F}_{\ell}$. $\square$

\subsection{A result of Roche on types}

We continue to assume $p$ is odd. We review a result of Roche, based on an explicit construction of types (in the sense of Bushnell and Kutzko) for principal series representations of split reductive groups. Thus, we let $I$ denote the Iwahori subgroup. That is, the inverse image of $B(\F)$ under the reduction map:
$$
I=\{\text{$g$ with $g$ (mod $\p$) $\in B(\F)$}\}\subset G_0(\mathcal{O}).
$$
Also, we introduce the normal subgroup $I_1$, which is the inverse image of $B_{u}(\F)$:
$$
I_1=\{\text{$g$ with $g$ (mod $\p$) $\in B_{u}(\F)$}\}\subset I.
$$
The quotient $I/I_1$ is naturally identified with $T(\F)$, by reducing the diagonal. If $\pi$ is an admissible representation of $G_0$, the finite abelian group $T(\F)$ acts on the $I_1$-invariants. Therefore this space has a decomposition into eigenspaces:
$$
\pi^{I_1}={\bigoplus}_{\chi_i,\chi} \pi^{I: \chi_1 \otimes \cdots \otimes \chi_n \otimes \chi}.
$$
Here the $\chi_i$ and $\chi$ run over all characters of $\F^*$, and the corresponding direct summand is the subspace of vectors transforming according to the associated character of $I$. We will need the following consequence of Roche's paper [Roc]:

\begin{thm}
Let $\pi$ be an irreducible admissible representation of $G_0(F)$. Then
$$
\pi^{I: \chi_1 \otimes \cdots \otimes \chi_n \otimes \chi}\neq 0
$$
if and only if $\pi$ is a subquotient of a principal series representation of the form
$$
\widetilde{\chi}_1 \times \cdots \times \widetilde{\chi}_n \rtimes \widetilde{\chi}
$$
for some quasi-characters $\widetilde{\chi}_i$ extending the $\chi_i$, and some $\widetilde{\chi}$ extending $\chi$.
\end{thm}

\noindent $Proof$. By inflation, we view $\chi_i$ and $\chi$ as characters of $\mathcal{O}^*$ trivial on the $1$-units $1+\p$. This defines an inertial equivalence class in the Bernstein spectrum,
$$
\frak{s}=[T, \widetilde{\chi}_1 \otimes \cdots \otimes \widetilde{\chi}_n \otimes \widetilde{\chi}].
$$
Here the $\widetilde{\chi}_i$ and $\widetilde{\chi}$ are arbitrary extensions of $\chi_i$ and $\chi$, respectively.
By the construction on page 367 in [Roc], we get an $\frak{s}$-type, according to Theorem 7.7 in [Roc]. This Theorem applies since $p$ is odd. By tameness, see Remark 4.2 in Roche's paper, the compact open subgroup of the $\frak{s}$-type is simply the Iwahori subgroup $I$. Moreover, the character of $I$ is the one coming from the $\chi_i$ and $\chi$ as above. By the very definition of a type, $\pi|_I$ contains this character if and only if $\frak{s}$ is the inertial class of $\pi$. That is, $\pi$ is a subquotient of a principal series
$$
\widetilde{\chi}_1 \times \cdots \times \widetilde{\chi}_n \rtimes \widetilde{\chi},
$$
where we may have to modify the characters by unramified twists. Done. $\square$

\medskip

\noindent It follows that $\pi^{I_1}$ is nonzero if and only if $\pi$ is a subquotient of a tamely ramified principal series (possibly unramified). As a continuation of Corollary 1, we get:

\begin{cor}
With $\pi$ as above, and the $\chi_i$ as in Corollary 1, we have for all $\chi$:
$$
\pi^{I: \chi_1 \otimes \cdots \otimes \chi_n \otimes \chi}\neq 0
$$
if and only if $\pi=\widetilde{\chi}_1 \times \cdots \times \widetilde{\chi}_n \rtimes \widetilde{\chi}$ for some
extensions $\widetilde{\chi}_i$ and $\widetilde{\chi}$ as above.
\end{cor}

\section{Modular forms and Hecke algebras}

For now, let us fix an arbitrary totally real number field $F$, of degree $d=[F:\Q]$.
We view $G_0$ as an algebraic $F$-group, and fix an inner form $G$ over $F$ such that
$$
G_{\infty}=G(F \otimes_{\Q}\R)=\prod_{v|\infty}G(F_v)
$$
is $compact$ modulo its center. The center $Z_G$ is isomorphic to $\G_m$ over $F$. The groups $G_0$ and $G$ are isomorphic over the algebraic closure $\bar{\Q}$, and throughout we fix an inner twisting $\psi$. By definition, this is an isomorphism over $\bar{\Q}$,
$$
\psi:G \rightarrow G_0,
$$
such that $\sigma \psi \circ \psi^{-1}$ is an inner automorphism of $G_0$ for all $\sigma \in \Gal(\bar{\Q}/F)$.
In this fashion, the equivalence class of $G$ corresponds to a cohomology class in
$$
H^1(F,G_0^{\text{ad}}).
$$
Here $G_0^{\text{ad}}$ is the adjoint group, and $H^1$ is non-abelian Galois cohomology. More concretely, the group $G$ is the unitary similitude group of an $n$-dimensional left $D$-vector space equipped with a non-degenerate hermitian form, where $D$ is a quaternion $F$-algebra [PR]. By choosing a basis, we can be even more concrete, and realize $G$ as an algebraic subgroup of $\GL_n(D)$: For any $F$-algebra $R$,
$$
G(R)=\{g \in \GL_n(D \otimes_{F}R): {^*g}\Phi g=\mu(g)\Phi\}.
$$
Here ${^*g}={^t\bar{g}}$, and $\Phi$ is non-degenerate and hermitian. That is, ${^*\Phi}=\Phi$. By an explicit calculation, it is not hard to show that the groups $G$ and $G_0$ are conjugate in $\GL_{2n}$ over $\bar{\Q}$.
We refer the reader to Lemma 2.3 in [Ghi] for the details. If $v$ is a non-archimedean place of $F$, there is a unique non-split inner form of $G_0$ over $F_v$ (the reduced norm on the division quaternion $F_v$-algebra is surjective, so all hermitian forms are equivalent). This inner form has rank $[n/2]$ modulo its center. If $v$ is archimedean, the hermitian forms are classified by their signature. Globally, there is always an $even$ number of places $v$ such that $G$ is non-split over $F_v$. At those places, the inner form can be prescribed.

\medskip

\noindent $Remark$. In the special case $n=2$, we may alternatively view $G$ as the similitude spin group of a totally definite quadratic form in five variables over $F$.

\subsection{Algebraic modular forms on $G$}

We wish to define modular forms on $G$ as in [Gro]. Thus, we first fix a central character $\omega$. This is naturally identified with a Hecke character of $F$, not necessarily unitary. Let $\mathcal{A}(\omega)$ be the space of $\omega$-central automorphic forms:
$$
\text{$f:G(F) \backslash G(\A_F)\rightarrow \C$, $\y$ $f(zg)=\omega(z)f(g)$, $\y$ $\forall z \in Z_G(\A_F)$.}
$$
These functions are assumed to be $G_{\infty}$-finite under right translation, and locally constant when restricted to the group of finite adeles $G(\A_{F}^f)$. In what follows, we simply say $f$ is smooth. The group $G(\A_F)$ acts on
the space $\mathcal{A}(\omega)$ by right translation, and it decomposes as a direct sum of irreducible automorphic representations $\pi=\pi_{\infty}\otimes \pi_f$ with $\omega_{\pi}=\omega$, each occurring with finite multiplicity.
We now fix an infinity type. That is, we fix an irreducible complex representation
$$
\text{$\xi:G_{\infty}\rightarrow \GL(V_{\xi})$, $\y$ $\xi=\otimes_{v|\infty}\xi_v$.}
$$
We assume it has central character $\omega_{\infty}$. This compatibility is necessary for the following spaces of modular forms to be non-trivial. First, let us introduce
$$
\mathcal{A}_{\xi}(\omega)=\Hom_{G_{\infty}}(\xi,\mathcal{A}(\omega))=
(\check{\xi}\otimes\mathcal{A}(\omega))^{G_{\infty}}.
$$
It carries an induced action of $G(\A_{F}^f)$, which is admissible, as we will see below.
The space $\mathcal{A}_{\xi}(\omega)$ is a direct sum of the $\pi_f$ such that $\xi\otimes \pi_f$ is automorphic, with central character $\omega$. Hence, the $\xi$-isotypic subspace of $\mathcal{A}(\omega)$ may be written as
$$
\mathcal{A}(\omega)[\xi]=\xi\otimes\mathcal{A}_{\xi}(\omega).
$$
We now want to realize $\mathcal{A}_{\xi}(\omega)$ as a space of vector-valued automorphic forms:

\begin{lem} As a $G(\A_{F}^f)$-representation, $\mathcal{A}_{\xi}(\omega)$ can be identified with the space of $\omega$-central $\check{V}_{\xi}$-valued smooth functions on the quotient $G(F)\backslash G(\A_F)$ such that
$$
\text{$\mathcal{F}(gg_{\infty}^{-1})=\check{\xi}(g_{\infty})\mathcal{F}(g)$, $\y$ $\forall g_{\infty}\in G_{\infty}$.}
$$
\end{lem}

\noindent $Proof$. To such an $\mathcal{F}$, we associate the map taking $u \in V_{\xi}$ to $u \circ \mathcal{F}$. $\square$

\medskip

\noindent By restriction to the finite adeles, we get: The space $\mathcal{A}_{\xi}(\omega)$ can be identified with the space of $\omega_f$-central $\check{V}_{\xi}$-valued smooth functions $\mathcal{F}$ on $G(\A_F^f)$ such that
$$
\text{$\mathcal{F}(\gamma_fg)=\check{\xi}(\gamma_{\infty})\mathcal{F}(g)$, $\y$ $\forall \gamma \in G(F)$.}
$$
Let $K$ be a compact open subgroup of $G(\A_F^f)$, and look at the $K$-invariants:
$$
\mathcal{A}_{\xi}(K,\omega)\simeq {\bigoplus}_{\pi:\pi_{\infty}=\xi}a_{\text{disc}}^G(\pi)\pi_f^K
$$
Such forms are said to have level $K$. Their union, as $K$ varies, is all of $\mathcal{A}_{\xi}(\omega)$. Moreover, the space above is finite-dimensional, thereby proving the aforementioned admissibility: Indeed, if we choose a set of representatives $g_1,\ldots,g_h$ for
$$
G(F)\backslash G(\A_F^f)/K,
$$
the map taking $\mathcal{F}$ to the $h$-tuple $(\mathcal{F}(g_1),\ldots,\mathcal{F}(g_h))$ defines an isomorphism
$$
\text{$\mathcal{A}_{\xi}(K,\omega)\simeq \bigoplus_{i=1}^h \check{V}_{\xi}^{\Gamma_i}$, $\y$ $\Gamma_i=G(F) \cap g_iG_{\infty}Kg_i^{-1}$.}
$$
Note that each $\Gamma_i$ is a finite subgroup of $G_{\infty}$. In fact, by shrinking $K$, one can often assume they are all trivial. For example, this is true if the projection of $K$ onto some $G_v$ does not contain any non-trivial elements of finite order. In this case, we say $K$ is sufficiently small. We will often make this assumption.

\medskip

\noindent $Remark$. We always assume $K$ is so small that $\omega_f$ is trivial on $K \cap Z_G(\A_F^f)$.

\subsection{Hecke algebras for $G$}

We let $\mathcal{H}(G)$ denote the Hecke algebra of $G(\A_F^f)$. That is, the smooth functions with compact support, equipped with the convolution product associated to some fixed Haar measure $dx$. It acts on $\mathcal{A}_{\xi}(\omega)$ in the usual way. That is,
$$
R(\phi)\mathcal{F}(g)=\int_{G(\A_F^f)}\phi(x)\mathcal{F}(gx)dx.
$$
Let $K$ be a compact open subgroup of $G(\A_F^f)$. If we normalize $dx$ such that $K$ has volume one, the characteristic function $e_K$ is an idempotent in $\mathcal{H}(G)$. Then
$$
\mathcal{H}(G,K)=e_K * \mathcal{H}(G) * e_K
$$
is the algebra of $K$-biinvariant compactly supported functions. Clearly it preserves
the level $K$ forms. If we fix a complete set of representatives $\{t\}$ for
$$
K\backslash G(\A_F^f)/K,
$$
the characteristic functions of the double cosets $KtK$ form a basis for $\mathcal{H}(G,K)$. In fact, the $\Z$-submodule they generate is a subring $\mathcal{H}_{\Z}(G,K)$. We introduce:
$$
\T_{\xi}(K,\omega)=\text{Im}\{R:\mathcal{H}_{\Z}(G,K)\rightarrow \text{End}_{\C}\mathcal{A}_{\xi}(K,\omega)\}.
$$
More generally, let $S$ be a finite set of finite places, and assume $K$ can be written as $K_SK^S$ where the factor $K^S$ is a product of hyperspecial subgroups $G(\mathcal{O}_v)$. Then the Hecke algebra $\mathcal{H}_{\Z}(G,K)$ factors correspondingly as a tensor product,
$$
\mathcal{H}_{\Z}(G,K)\simeq \mathcal{H}_{\Z}(G_S,K_S)\otimes\mathcal{H}_{\Z}^S(G,K),
$$
and we let $\T_{\xi}^S(K,\omega)$ denote the image of the prime-to-$S$ part $\mathcal{H}_{\Z}^S(G,K)$ in the above endomorphism algebra. This factor of the Hecke algebra is commutative.

\subsection{Arithmetic models}

Next, we define models of the spaces $\mathcal{A}_{\xi}(K,\omega)$ over number fields and finite fields. Thus, we assume our inner twisting $\psi$ is defined over the number field $L$ inside $\bar{\Q}$. By enlarging $L$, we may assume $\omega_f$ maps into $L^*$, and that $L$ contains $F$. The complexification of $\xi$ is an algebraic representation of the complex group
$$
G(F \otimes_{\Q}\C)=\prod_{v|\infty}G(F_v \otimes_{\R}\C).
$$
Such an irreducible representation $\otimes_{v|\infty}\xi_v$ is characterized by its highest weight. This is a dominant character of $T$. In our case, since $G$ is isomorphic to $G_0$ over $\C$, the highest weight corresponds to a tuple $a=(a_v)_{v|\infty}$ where for each $v$:
$$
\text{$a_v=(a_{v,1},\ldots,a_{v,n},\widetilde{a}_v)\in \Z^{n+1}$, $\y$ $a_{v,1}\geq \cdots \geq a_{v,n}\geq 0$.}
$$
Again, by enlarging $L$, we may assume this algebraic representation is defined over $L$, viewed as a subfield of $\C$ via our fixed embedding $i_{\infty}$. We denote it by
$$
\text{$\xi: \widetilde{G}_{/L}\rightarrow \GL (V_{\xi}(L))$, $\y$ $\widetilde{G}=\text{Res}_{F/\Q}(G)$.}
$$
In particular, the group $G(F)$ acts naturally on $V_{\xi}(L)$: It is a subgroup of $\widetilde{G}(L)$ via the diagonal embedding. This allows us to introduce an $L$-rational structure on the spaces of modular forms. Indeed, we make the following definition:

\begin{df}
We define the $L$-rational modular forms $\mathcal{A}_{\xi}(K,\omega,L)$ to be the space of $\omega_f$-central $\check{V}_{\xi}(L)$-valued $K$-invariant functions $\mathcal{F}$ on $G(\A_F^f)$ such that
$$
\text{$\mathcal{F}(\gamma_fg)=\check{\xi}(\gamma)\mathcal{F}(g)$, $\y$ $\forall \gamma \in G(F)$.}
$$
\end{df}

\noindent It carries an action of the Hecke algebra $\T_{\xi}(K,\omega)$. We now complete the above space $\ell$-adically, for a fixed prime $\ell$. Thus, let $\lambda$ be the place of $L$ determined by the embedding $i_{\ell}$, and let $L_{\lambda}$ be the corresponding completion of $L$. We let
$$
\mathcal{A}_{\xi}(K,\omega,L_{\lambda})=L_{\lambda}\otimes_L \mathcal{A}_{\xi}(K,\omega,L).
$$
Of course, this can be identified with smooth $\check{V}_{\xi}(L_{\lambda})$-valued functions on $G(\A_F^f)$ as above. However, for the purpose of defining integral models and modular forms mod $\ell$, it is better to work with a different model. Namely, the following:

\begin{lem}
As a $\T_{\xi}(K,\omega)$-module, the space $\mathcal{A}_{\xi}(K,\omega,L_{\lambda})$ can be identified with the space of $\check{V}_{\xi}(L_{\lambda})$-valued functions $\mathcal{F}$ on the quotient $G(F)\backslash G(\A_F^f)$ such that
$$
\text{$\mathcal{F}(gk)=\check{\xi}(k_{\lambda}^{-1})\mathcal{F}(g)$, $\y$ $\forall k \in K$,}
$$
and such that the prime-to-$\lambda$ part of $Z_G(\A_F^{f})$ acts on $\mathcal{F}$ by the character $\omega_f^{\lambda}$.
\end{lem}

\noindent $Remarks$. We abuse notation, and use $\lambda$ also to denote the place of $F$ defined by $i_{\ell}$.
The Hecke action in this model, is given by the following explicit formula:
$$
R(e_{KtK})\mathcal{F}(g)={\sum}_{x \in KtK/K}\check{\xi}(x_{\lambda})\mathcal{F}(gx).
$$
\medskip
\noindent $Proof$. For a function $\mathcal{F}$ in the former model, the function $g \mapsto \check{\xi}(g_{\lambda}^{-1})\mathcal{F}(g)$ is left-invariant under $G(F)$, and is easily seen to belong to the latter model - and vice versa: The compatibility between $\omega$ and the central character of $\xi$ ensures that  $F_{\lambda}^*$ acts on the form $\mathcal{F}$ by the central character $\omega_{\lambda}$ as wanted. $\square$

\medskip

\noindent Now, to define integral models, over the valuation ring $\mathcal{O}_{\lambda}$, we need a lattice
preserved by all $k_{\lambda}$ for $k \in K$. We denote the group of such elements by $K_{\lambda}$. For this, we assume $\widetilde{G}$ extends to a smooth affine group scheme over $\mathcal{O}_{\lambda}$, and
$$
\xi:\widetilde{G}_{/\mathcal{O}_{\lambda}}\rightarrow \GL(V_{\xi}(\mathcal{O}_{\lambda}))
$$
is an extension of the previous $\xi$. We assume $\widetilde{G}(\mathcal{O}_{\lambda})$ contains $K_{\lambda}$.
Consequently,
$$
\mathcal{A}_{\xi}(K,\omega,L_{\lambda})=L_{\lambda}\otimes_{\mathcal{O}_{\lambda}} \mathcal{A}_{\xi}(K,\omega,\mathcal{O}_{\lambda}),
$$
where $\mathcal{A}_{\xi}(K,\omega,\mathcal{O}_{\lambda})$ is the lattice of $\check{V}_{\xi}(\mathcal{O}_{\lambda})$-valued functions $\mathcal{F}$ as above. We immediately observe that this lattice is not preserved the full Hecke algebra $\T_{\xi}(K,\omega)$, unless $\xi$ is trivial. Only the prime-to-$\lambda$ part $\T_{\xi}^{\lambda}(K,\omega)$ preserves integrality. Finally, to define modular forms mod $\ell$, we tensor by $\F_{\lambda}$. That is,
$$
\mathcal{A}_{\xi}(K,\omega,\F_{\lambda})=\F_{\lambda}\otimes_{\mathcal{O}_{\lambda}} \mathcal{A}_{\xi}(K,\omega,\mathcal{O}_{\lambda}).
$$
This can of course be canonically identified with the space of $\check{V}_{\xi}(\F_{\lambda})$-valued functions $\mathcal{F}$ satisfying the transformation properties in the previous lemma.

\section{Congruences}

In this section, we will prove Theorem A in the introduction. Let us briefly review the setup: Start with an automorphic representation $\pi$ of $G(\A_F)$ with
$$
\text{$\omega_{\pi}=\omega$, $\y$ $\pi_{\infty}=\xi$, $\y$ $\pi_f^I \neq 0$,}
$$
where $I$ is a compact open subgroup of $G(\A_F^f)$ of the form $I_wK^w$ for a fixed finite place $w$ of $F$ where $G$ splits. Here $I_w$ is an Iwahori subgroup of $G(F_w)$, and $K^w$ is some compact open subgroup prime-to-$w$. Moreover, we assume:
$$
\text{$N(\p_w)\equiv 1$ (mod $\ell$), $\y$ $\ell>2n$.}
$$
Let us furthermore assume $K^w$ is factorizable of the form $\prod_{v\neq w}K_v$, where $K_v$ is hyperspecial for $v \notin S_0$. We let $S_{\ell}$ denote the places of $F$ above $\ell$, and consider
$$
S=\{w\}\cup S_0\cup S_{\ell}.
$$
Then the Hecke algebra $\T_{\xi}^S(I,\omega)$ acts on $\pi_f^I$ via an algebra character, denoted
$$
\eta_{\pi}: \T_{\xi}^S(I,\omega) \rightarrow \C.
$$
It maps into $\bar{\Q}$, since the Hecke algebra preserves $\mathcal{A}_{\xi}(I,\omega,L)$. In fact, via $i_{\ell}$, it maps into $\bar{\Z}_{\ell}$. Indeed, the previously defined integral structures are preserved. It therefore makes sense to reduce $\eta_{\pi}$ mod $\ell$. We denote this $\bar{\F}_{\ell}$-valued character by $\bar{\eta}_{\pi}$, and let
$\frak{m}$ be its kernel. This is a maximal ideal of the algebra $\T_{\xi}^S(I,\omega)$.

\medskip

\noindent $Goal$: Find an automorphic representation $\widetilde{\pi}\equiv \pi$ (mod $\ell$) of $G(\A_F)$ such that
$$
\text{$\omega_{\widetilde{\pi}}=\omega$, $\y$ $\widetilde{\pi}_{\infty}=\xi$, $\y$ $\widetilde{\pi}_f^{K^w} \neq 0$,}
$$
and such that $\widetilde{\pi}_w$ is a tamely ramified principal series representation of $G(F_w)$.

\medskip

\noindent Here, the  notation $\widetilde{\pi}\equiv \pi$ (mod $\ell$) means that $\bar{\eta}_{\widetilde{\pi}}=\bar{\eta}_{\pi}$, where the character
$$
\eta_{\widetilde{\pi}}: \T_{\xi}^S(I,\omega) \rightarrow \C
$$
gives the action of $\T_{\xi}^S(I,\omega)$ on the $K^w$-invariants of $\widetilde{\pi}$. To reach the goal, we use the characters $\chi_i$ of the unit group $\mathcal{O}_w^*$ from Corollary 1. Moreover, fix an arbitrary $\chi\equiv 1$ (mod $\ell$). Having fixed the $\{\chi_i,\chi\}$, we introduce the subspace
$$
\mathcal{A}_{\xi}^{\{\chi_i,\chi\}}(I,\omega)
$$
of forms $\mathcal{F}$ invariant under $I_{1w}K^w$, on which $I_w$ acts by $\otimes\chi_i\otimes \chi$. We need:
$$
\mathcal{A}_{\xi}^{\{\chi_i,\chi\}}(I,\omega)_{\frak{m}}\neq 0.
$$
Indeed, if this is the case, by Corollary 2 there is a $\widetilde{\pi}$ with $\bar{\eta}_{\widetilde{\pi}}=\bar{\eta}_{\pi}$ such that
$$
\widetilde{\pi}_w=\widetilde{\chi}_1 \times \cdots \times \widetilde{\chi}_n \rtimes \widetilde{\chi}
$$
for some extensions $\widetilde{\chi}_i$ and $\widetilde{\chi}$. It suffices to show non-vanishing at $\frak{m}$ mod $\ell$:
$$
\mathcal{A}_{\xi}^{\{\chi_i,\chi\}}(I,\omega,\F_{\lambda})_{\frak{m}}\simeq\F_{\lambda}\otimes_{\mathcal{O}_{\lambda}}
\mathcal{A}_{\xi}^{\{\chi_i,\chi\}}(I,\omega,\mathcal{O}_{\lambda})_{\frak{m}}.
$$
Now recall from Corollary 1 that all the $\chi_i$ are trivial mod $\ell$. Therefore,
$$
\mathcal{A}_{\xi}^{\{\chi_i,\chi\}}(I,\omega,\F_{\lambda})_{\frak{m}}=
\mathcal{A}_{\xi}^{\{1,1\}}(I,\omega,\F_{\lambda})_{\frak{m}},
$$
and this is simply the space of $I$-invariants. Here $\frak{m}$ does occur. We are done.

\section{Jacquet-Langlands transfer for $n=2$}

To get congruences between automorphic forms on $G_0$, we need to establish certain cases of functoriality, analogous to the Jacquet-Langlands correspondence for $\GL_2$ and its inner forms. For this, we need to impose the two conditions:
$$
\text{$n=2$, $\y$ $d=[F:\Q]$ is even.}
$$
The first condition is dictated by the present state of the trace formula for symplectic groups. The second condition guarantees the existence of a totally definite quaternion algebra $D$ over $F$, which is split at all the finite places:
$$
\text{$D\otimes_{\Q}\R\simeq\HH^d$, $\y$ $D_v\simeq M_2(F_v)$,}
$$
for all finite $v$. Here $\HH$ denotes the Hamilton quaternions over $\R$. These conditions determine $D$ uniquely up to isomorphism. We take the hermitian form $\Phi$ to be the identity, and look at the unitary similitude group: For $F$-algebras $R$,
$$
G(R)=\{g \in \GL_2(D \otimes_F R): {^*g}g=\mu(g)I\}.
$$
Then $G_{\infty}$ is compact modulo its center, and $G$ is split over $F_v$ for all finite $v$.

\subsection{The trace formula for $G$}

The group $G$ is anisotropic modulo its center $Z_G$, so the trace formula takes its simplest form. We keep the central character $\omega$ fixed, and look at the space
$$
C_c^{\infty}(G(\A_F),\check{\omega})
$$
of smooth $\check{\omega}$-central test functions $\phi$ on $G(\A_F)$, assumed to be compactly supported modulo the center. Each such function $\phi$ defines a compact operator
$$
R(\phi)f(g)=\int_{G^{\text{ad}}(\A_F)}\phi(x)f(gx)dx
$$
on $\mathcal{A}(\omega)$. It turns out to be of trace class, and the trace formula computes its trace in two different ways. First, one obviously has the spectral expansion:
$$
I_{\text{disc}}^G(\phi)\overset{\text{df}}{=}\tr R(\phi)={\sum}_{\pi: \omega_{\pi}=\omega}a_{\text{disc}}^G(\pi)\tr \pi(\phi),
$$
in terms of the automorphic spectrum of $G$. One compares it with the geometric expansion of the trace:
First, observe that the operator $R(\phi)$ is given by a kernel,
$$
R(\phi)f(g)=\int_{Z_G(\A_F)G(F)\backslash G(\A_F)}K_{\phi}(g,x)f(x)dx,
$$
where we have introduced
$$
K_{\phi}(g,x)={\sum}_{\gamma \in G^{\text{ad}}(F)}\phi(g^{-1}\gamma x).
$$
To get the trace, we integrate $K_{\phi}$ over the diagonal. A short calculation shows:
$$
\tr R(\phi)=\sum_{\{\gamma\}}\text{vol}(Z_G(\A_F)G_{\gamma}(F)\backslash G_{\gamma}(\A_F))\cdot O_{\gamma}(\phi).
$$
Here $\{\gamma\}$ is a set of representatives for the conjugacy classes in $G^{\text{ad}}(F)$, and $G_{\gamma}$ is the centralizer of $\gamma$ in $G$. Moreover, $O_{\gamma}$ denotes the orbital integral. That is,
$$
O_{\gamma}(\phi)=\int_{G_{\gamma}(\A_F)\backslash G(\A_F)}\phi(x^{-1}\gamma x)dx.
$$

\subsection{Stable orbital integrals}

One would like to rewrite the geometric side in terms of $stable$ orbital integrals. To define these, let us first recall the notion of stable conjugacy: Eventually we want to compare the geometric side of the trace formula for $G$ with that of the (much more complicated) trace formula for $G_0$, in order to compare automorphic spectra. However, the two groups are only isomorphic over $\bar{\Q}$, so one can only compare conjugacy classes over $\bar{\Q}$. We work locally at a fixed place $v$ of $F$:

\begin{df}
Elements $\gamma$ and $\gamma'$ in $G(F_v)$ are $stably$ conjugate if and only if
$$
\text{$\exists$ $g \in G(\bar{F}_v)$: $\gamma'=g^{-1}\gamma g$ (mod center).}
$$
\end{df}

\noindent In the case of $\GL_n$, this turns out to coincide with usual conjugacy. However, for $G_0$ it is coarser. We have used the fact that the derived group $G^{\text{der}}$ is simply connected; otherwise, the definition has to be altered slightly. A stable conjugacy class contains only finitely many conjugacy classes, parameterized by
$$
\ker\{H^1(F_v,G_{\gamma})\rightarrow H^1(F_v,G)\}.
$$
Again, $H^1$ is non-abelian Galois cohomology. The parametrization is not difficult; one looks at the cohomology class of the cocycle $\sigma(g)g^{-1}$. In general, for arbitrary semisimple $\gamma$, the centralizer $G_{\gamma}$ is a connected reductive group. If $\gamma$ is regular, $G_{\gamma}$ is a torus. It is connected, again since
$G^{\text{der}}$ is simply connected. One can show that the above kernel is a finite abelian subgroup.

\begin{df}
Suppose $\gamma'$ is stably conjugate to $\gamma$. Then $G_{\gamma'}$ is the inner form of $G_{\gamma}$ corresponding to the class in $H^1(F_v,G_{\gamma}^{\text{ad}})$ given by the conjugacy class of $\gamma'$ in the above parametrization. Let $e(G_{\gamma'})$ be the Kottwitz sign [Kot]. Then
$$
SO_{\gamma}(\phi)=\sum_{\gamma'}e(G_{\gamma'})O_{\gamma'}(\phi)
$$
is the stable orbital integral. Here $\gamma'\sim\gamma$ runs over the conjugacy classes.
\end{df}

\noindent Globally, stable orbital integrals are defined as products of local ones. More generally, one defines $\kappa$-orbital integrals for certain characters $\kappa$.The stable case corresponds to $\kappa=1$. A distribution $I$ on $G(\A_F)$ is called $stable$ if it is supported on the stable orbital integrals. That is, for all test functions $\phi$ as above,
$$
\text{$SO_{\gamma}(\phi)=0$, $\forall$ semisimple $\gamma \in G(F)$ $\Longrightarrow$ $I(\phi)=0$.}
$$

\subsection{Elliptic endoscopic triples for $G$}

Unfortunately, one cannot express the geometric side of the trace formula for $G$ purely in terms of stable orbital integrals on $G$. The distribution $I_{\text{disc}}^G$ turns out to be unstable. The solution to this problem is to write $I_{\text{disc}}^G$ as a sum of stable distributions on the so-called elliptic endoscopic groups. Actually, to be precise, one looks at equivalence classes of triples $(H,s,\zeta)$, consisting of data:
\begin{itemize}
\item $H$ is a quasi-split connected reductive group over $F$,\\
\item $s$ is a semisimple element in the dual group $\hat{G}$,\\
\item $\zeta: \hat{H}\hookrightarrow \hat{G}$ is a homomorphism with image $Z_{\hat{G}}(s)^0$.
\end{itemize}

\noindent Such a triple must satisfy two conditions that we will not state precisely, see page 20 in [Har]. First, $\zeta$ must be Galois-equivariant up to $\hat{G}$-conjugation. Secondly, $s$ must satisfy some local triviality condition. Finally, the triple is $elliptic$ if
$$
\text{$\zeta(Z(\hat{H})^{\Gamma}) Z(\hat{G})/Z(\hat{G})$, $\y$ $\Gamma=\Gal(\bar{\Q}/F)$,}
$$
is finite. There is a natural equivalence relation on such triples, see Definition 5.2 in [Har]. In our case, $G$ admits precisely two elliptic endoscopic triples up to equivalence. The first one being the quasi-split inner form $G_0$ equipped with $s=1$ and $\zeta=\text{Id}$. The connected reductive group of the second triple is
$$
H=[\GL(2)\times \GL(2)]/\G_m
$$
Here $\G_m$ is a subgroup of the center via the embedding $z \mapsto (z,z^{-1})$. Moreover,
$$
s=\begin{pmatrix}1 & & & \\ & -1 & & \\ & & -1 & \\ & & & 1\end{pmatrix} \in \hat{G}=\GSp_4(\C).
$$
The dual group $\hat{H}$ is the subgroup of $\GL_2(\C)\times \GL_2(\C)$ consisting of pairs of matrices with the same determinant. We may take the homomorphism $\zeta$ to be
$$
\zeta:
\begin{pmatrix}a & b \\ c & d\end{pmatrix}\times
\begin{pmatrix}\widetilde{a} & \widetilde{b} \\ \widetilde{c} & \widetilde{d}\end{pmatrix}\mapsto \begin{pmatrix}\widetilde{a} & & & \widetilde{b}\\ & a & b & \\ &c & d & \\ \widetilde{c} & & & \widetilde{d}\end{pmatrix}.
$$
It is straightforward to check that its image is indeed the centralizer of $s$ in $\hat{G}$.

\subsection{Transfer and the fundamental lemma}

In this section, in the beginning, we fix a place $v$ and let $G=G_v$ and so on. We wish to make precise, how one transfers a given test function $\phi$ on $G$, or $G_0$, to families of functions on the endoscopic groups $G_0$ and $H$. Let us first explain how to go from $G$ to $G_0$: The set of (semisimple) stable conjugacy classes in $G$ embeds into the set of stable conjugacy classes in $G_0$. Namely, let $\gamma \in G$ be semisimple, and look at the $G_0$-conjugacy class of $\psi(\gamma)$ over $\bar{F}_v$. Recall that we have fixed an inner twisting $\psi$. By a theorem of Kottwitz and Steinberg, since $G_0$ is quasi-split and $G_0^{\text{der}}$ is simply connected, the
conjugacy class of $\psi(\gamma)$ intersects the $F_v$-rational points $G_0$ in a stable conjugacy class. That is,
$$
\text{$\exists$ $\gamma_0 \in G_0$ such that $\gamma_0 \sim_{G_0(\bar{F}_v)} \psi(\gamma)$.}
$$
We write $\gamma \leftrightarrow \gamma_0$ in this case, and say that $\gamma_0$ is an $image$ of $\gamma$. For a given $\gamma_0$, when such a $\gamma$ exists, we say that $\gamma_0$ $occurs$ in $G$. We need the following result:

\begin{thm}
For every $\phi \in C_c^{\infty}(G,\check{\omega}_v)$ there exists a $\phi^{G_0} \in C_c^{\infty}(G_0,\check{\omega}_v)$ with
$$
SO_{\gamma_0}(\phi^{G_0})=SO_{\gamma}(\phi)
$$
when $\gamma \leftrightarrow \gamma_0$, and such that $SO_{\gamma_0}(\phi^{G_0})=0$ whenever $\gamma_0$ does not occur in $G$.
\end{thm}

\noindent $Proof$. In the archimedean case, this was proved by Shelstad in [Sh1]. In the non-archimedean case, it is a deep Theorem of Waldspurger, proved in [Wal]. $\square$

\medskip

\noindent In this situation, we say that $\phi$ and $\phi^{G_0}$ have $matching$ orbital integrals. Of course, $\phi^{G_0}$ is not unique, but $I(\phi^{G_0})$ is well-defined for stable distributions $I$.
Next, let us explain how to go from $G$ to $H$. The discussion with $G$ replaced by $G_0$ is similar, but simpler.
As above, take a semisimple element $\gamma \in G$. We first define what it means for a semisimple $\gamma_H \in H$ to be an $image$ of $\gamma$: Choose a maximal torus $T_H$ in the centralizer $H_{\gamma_H}$, containing $\gamma_H$. There is a canonical $G$-conjugacy class of embeddings of $T_H$ into $G$ defined over $\bar{F}_v$. Choose one, say
$$
j:T_H \hookrightarrow G,
$$
and assume it maps into the maximal torus $T$ in $G$. In analogy with the above discussion, we consider the stable conjugacy class of $j(\gamma_H)$. If it contains $\gamma$, we say that $\gamma_H$ is an image of $\gamma$. Alternative terminology; say $\gamma$ $comes$ from $\gamma_H$.

\begin{df}
The roots $\Phi$ of $G$ relative to $T$ define characters of $T_H$ by pulling back via $j$. We say that a semisimple element $\gamma_H \in H$ is $(G,H)$-regular if $\alpha(\gamma_H)\neq 1$ for all $\alpha \in \Phi$ such that the $j$-restriction $\alpha_{T_H}$ is not a root of $H$.
\end{df}

\noindent The following is the Langlands-Shelstad transfer conjecture, known in our case:

\begin{thm}
For every $\phi \in C_c^{\infty}(G,\check{\omega}_v)$ there exists a $\phi^{H} \in C_c^{\infty}(H,\check{\omega}_v)$ with
$$
SO_{\gamma_H}(\phi^H)=\sum_{\gamma}\Delta_{G,H}(\gamma_H,\gamma)e(G_{\gamma})O_{\gamma}(\phi)
$$
for all $(G,H)$-regular $\gamma_H$. The sum is over $\gamma$ coming from $\gamma_H$, up to conjugacy.
\end{thm}

\noindent $Proof$. Again, the archimedean case was worked out completely for all groups by Shelstad in [Sh1] and [Sh2]. The non-archimedean case follows from work of Hales, [Ha1] and [Ha2], combined with either [LS] or [Wal]. $\square$

\medskip

\noindent $Remarks$. Once we have normalized the Langlands-Shelstad transfer factors $\Delta_{G,H}$, compatible measures must be used on both sides. Fortunately, we will not need a precise definition of $\Delta_{G,H}$, so we will not give it here. The reader is referred to Chapter 7 in [Har] for a nice review of the basic properties. We also note that, since $H$ has no endoscopy, the left-hand side of the transfer identity is simply the standard orbital integral $O_{\gamma_H}(\phi^H)$, up to a sign at most.

\medskip

\noindent Globally, the test functions are spanned by pure tensors $\phi=\otimes \phi_v$, where $\phi_v$ is the characteristic function of a hyperspecial subgroup $K_v$ for almost all $v$; they are all conjugate. We would like to transfer $\phi$ to a function $\phi^H$ by the formula
$$
\phi^H\overset{\text{df}}{=}\otimes_v \phi_v^{H_v}.
$$
However, for this to make sense, we need to know that we can take $\phi_v^{H_v}$ to be the characteristic function of a hyperspecial subgroup $K_v^H$ for almost all $v$. This is the so-called standard fundamental lemma, proved by Hales in the case $n=2$:

\begin{thm}
$e_{K_v}^{H_v}=e_{K_v^H}$ (up to a nonzero constant depending on the measures).
\end{thm}

\noindent $Proof$. This is the main result of [Ha1]. $\square$

\medskip

\noindent $Remark$. We mention an amazing observation of Waldspurger: For general groups, the fundamental lemma in fact implies the transfer conjecture [Wal].

\medskip

\noindent To recapitulate, we can transfer a $global$ function $\phi$ on $G(\A_F)$ to a function $\phi^H$ on $H(\A_F)$, locally with matching orbital integrals, which is well-defined on stable distributions. Similarly for functions on the group $G_0(\A_F)$. Furthermore,
$$
\phi^{G_0}\overset{\text{df}}{=}\otimes_v \phi_v^{G_{0,v}}
$$
defines the transfer of $\phi$ to $G_0(\A_F)$, well-defined since $G=G_0$ at almost all $v$.

\subsection{Character relations for the endoscopic lift}

We go back to the local notation from the previous section. Thus, we fix a place $v$ and let
$G=G_v$ and so on. In this section, we discuss the local functorial lift pertaining to the $L$-homomorphism extending the map $\zeta$: The lift takes an irreducible admissible representation $\rho$ of $H$, and associates to it an $L$-packet of representations of $G_0$. If a certain relevancy condition is satisfied, this packet descends to $G$. The representation $\rho$ is of the form $\rho_1\otimes \rho_2$, where the $\rho_i$ are representations of $\GL_2$ with the same central character. The distribution $\tr\rho$ is stable, since $H$ has no endoscopy. By general results of Arthur [Ar2] in the non-archimedean case, and of Shelstad [Sh1] in the archimedean case,
$$
\text{$\tr \rho(\phi^H) =\sum_{\pi}\Delta_{G,H}(\rho,\pi)\tr\pi(\phi)$, $\y$ $\forall$ $\phi \in C_c^{\infty}(G,\check{\omega}_v)$,}
$$
for some coefficients $\Delta_{G,H}(\rho,\pi)$ one should think of as spectral analogues of the transfer factors. It is zero for almost all irreducible admissible representations $\pi$ of $G$. Similar expansions hold when $G$ is replaced by $G_0$. In the non-archimedean case, the above expansion has been made very precise by Weissauer in the two preprints [We1] and [We2]: For example, when $\rho$ is in the $discrete$ series,
$$
\tr\rho(\phi^H)=\tr\pi_+(\rho)(\phi)-\tr\pi_{-}(\rho)(\phi).
$$
Here $\{\pi_{\pm}(\rho)\}$ is the associated $L$-packet of representations of $G$ (recall that $G$ and $G_0$ are isomorphic at $all$ finite places). The representation $\pi_{+}(\rho)$ is always generic, in the sense that it has a Whittaker model, whereas the other member $\pi_{-}(\rho)$ is always non-generic. One has the following qualitative description:
\begin{itemize}
\item $\rho_1\neq \rho_2$ ($\rho$ $regular$): $\pi_+(\rho)$ is a discrete series, $\pi_-(\rho)$ is supercuspidal.\\
\item $\rho_1=\rho_2$ ($\rho$ $invariant$): $\{\pi_{\pm}(\rho)\}$ are limits of discrete series ($\Rightarrow$ tempered).
\end{itemize}
We will not need it in this paper, but when the residue characteristic of $F_v$ is $odd$, the $\pi_{\pm}(\rho)$ have a nice description in terms of $\theta$-correspondence for similitude groups. This comes about, by identifying $H$ and its inner form $\breve{H}$ with $\GSO$'s:
$$
\text{$H\simeq \GSO(V_{\text{split}})$, $\y$ $\breve{H}=[D^*\times D^*]/\G_m \simeq \GSO(V_{\text{aniso}})$.}
$$
Here $V_{\text{split}}$ and $V_{\text{aniso}}$ are the split and anistropic quaternary quadratic spaces of discriminant one, respectively. Since $\rho$ is in the discrete series, it corresponds to a unique representation $\breve{\rho}$ of $\breve{H}$ under the Jacquet-Langlands correspondence for $\GL_2$. By [Rob], both $\rho$ and $\breve{\rho}$ have unique extensions to the corresponding $\GO$'s participating in the $\theta$-correspondence with $\GSp_4$. We denote them by $\rho^+$ and $\breve{\rho}^+$. In the regular case, they are simply the induced representations. Then:
$$
\text{$\pi_+(\rho)=\theta(\rho^+)$, $\y$ $\pi_-(\rho)=\theta(\breve{\rho}^+)$.}
$$
For this, see Proposition 1 on page 14 in [W2]. In the case of $even$ residue characteristic, one can can still get an interpretation in terms of $\theta$-correspondence by looking at the $globally$ $relevant$ subset. Finally, when $\rho$ is a $principal$ series,
$$
\tr\rho(\phi^H)=\tr\pi_+(\rho)(\phi),
$$
for a generic representation $\pi_+(\rho)$ (tempered if $\rho$ is tempered). When one of the $\rho_i$ are one-dimensional, the support of the expansion of $\tr\rho$ does not define $L$-packets, but $A$-packets. Their description is given explicitly in [W1]. One should mention that, the above results of Weissauer are complete analogues of the results known in the archimedean case, due to the work of Shelstad [Sh1].
Globally, one immediately deduces exact analogues of the above character expansions: If $\rho=\otimes\rho_v$ is an automorphic representation of $H$, then we have
$$
\text{$\tr \rho(\phi^H) =\sum_{\pi}\Delta_{G,H}(\rho,\pi)\tr\pi(\phi)$, $\y$ $\forall$ $\phi \in C_c^{\infty}(G(\A_F),\check{\omega})$,}
$$
where $\Delta_{G,H}(\rho,\pi)$ is the product of all the local coefficients $\Delta_{G_v,H_v}(\rho_v,\pi_v)$.
Here $\pi$ runs over the global $L$-packet for $G$ associated with $\rho$, that is, the tensor product of all the local $L$-packets. It is possibly empty for some $\rho$.

\subsection{Stabilization of the trace formula}

A key observation is that $H$ itself has $no$ endoscopy; stable conjugacy is the same as conjugacy. Indeed, $H$ is essentially just a product of two copies of $\GL_2$. Nevertheless, the trace formula for $H$ is more complicated than that for $G$, due to the existence of parabolic subgroups defined over $F$. One has to build the continuous spectrum from Eisenstein series, and subtract this contribution to get the trace of $R(\phi)$ on the discrete spectrum. In the case of $\GL_2$ this is done in detail in [GJ], and Arthur dealt with arbitrary rank one groups in [Ar1]. In addition to the trace of $R(\phi)$ on the discrete spectrum, there are certain so-called scattering terms occurring discretely on the spectral side of the trace formula. We collect all these terms, and denote the resulting distribution by $I_{\text{disc}}^H$. From what we have just said, this is a stable distribution. It can be written as follows:
$$
\text{$I_{\text{disc}}^H(\phi)={\sum}_{\rho:\omega_{\rho}=\omega}a_{\text{disc}}^H(\rho)\tr \rho(\phi)$, $\y$ $\forall$ $\phi \in C_c^{\infty}(H(\A_F),\check{\omega})$.}
$$
Here $\rho$ runs over all discrete automorphic representations of $H$ with $\omega_{\rho}=\omega$. It is of the form $\rho_1\otimes \rho_2$, where the $\rho_i$ are automorphic representations of $\GL_2$ with central character $\omega$.
In general, the complex coefficient $a_{\text{disc}}^H(\rho)$ may not be equal to the multiplicity of $\rho$. However, in the generic case where no $\rho_i$ is a character, they $are$ equal. The trace formula for $G_0$ is even more complicated. There are now three proper parabolic subgroups, up to conjugacy, each contributing to the continuous spectrum. Again, besides the trace of $R(\phi)$ on the discrete spectrum, we get additional discretely occurring terms on the spectral side. They are what Arthur refers to as surviving remnants of Eisenstein series. As before, we construct a distribution $I_{\text{disc}}^{G_0}$ out of these terms, and write it as:
$$
\text{$I_{\text{disc}}^{G_0}(\phi)={\sum}_{\Pi:\omega_{\Pi}=\omega}a_{\text{disc}}^{G_0}(\Pi)\tr \Pi(\phi)$, $\y$ $\forall$ $\phi \in C_c^{\infty}(G_0(\A_F),\check{\omega})$.}
$$
Here $\Pi$ varies over all discrete automorphic representations of $G_0$ with central character $\omega$. The coefficient $a_{\text{disc}}^{G_0}(\Pi)$ equals the multiplicity of $\Pi$, if $\Pi$ is cuspidal but not CAP. Recall that a cuspidal automorphic representation $\Pi$ is CAP with respect to a parabolic $P$, if it is nearly equivalent to the constituents of a representation induced from a cusp form on a Levi factor $M$ of $P$. For $\GL_n$, this notion is vacuous by a result of Shalika. On the other hand, for $G_0$, the CAP forms have been constructed via theta correspondence by Piatetski-Shapiro [PS] and Soudry [Sou]. As with $I_{\text{disc}}^G$, the distribution $I_{\text{disc}}^{G_0}$ is also unstable. However, if we subtract suitable endoscopic error terms, we can $make$ it stable:

\begin{thm}
The distribution $S_{\text{disc}}^{G_0}(\phi)\overset{\text{df}}{=}I_{\text{disc}}^{G_0}(\phi)-\frac{1}{4}I_{\text{disc}}^H(\phi^H)$ is stable.
\end{thm}

\noindent $Proof$. This is a special case of a Theorem of Arthur announced in [Ar3]. $\square$

\medskip

\noindent We thus arrive at the stable trace formula for $G$:

\begin{thm}
$I_{\text{disc}}^G(\phi)=S_{\text{disc}}^{G_0}(\phi^{G_0})+\frac{1}{4}I_{\text{disc}}^H(\phi^H)$.
\end{thm}

\noindent $Proof$. This was first proved by Kottwitz [Kot] and Langlands [La2], but can be viewed also as a special case of the aforementioned work of Arthur [Ar3]. $\square$

\subsection{A trace identity at infinity}

Finally, we are ready to compare the automorphic spectra of $G$ and $G_0$. The link is the stable distribution
$S_{\text{disc}}^{G_0}$. We insert the character expansions of the $\tr\rho$ into $I_{\text{disc}}^H$, in order to express the contribution from $H$ to the stable trace formula in terms of the spectrum of $G$. For all test functions $\phi$ on $G$, we get:
$$
S_{\text{disc}}^{G_0}(\phi^{G_0})=I_{\text{disc}}^G(\phi)-\frac{1}{4}I_{\text{disc}}^H(\phi^H)={\sum}_{\pi: \omega_{\pi}=\omega}\widetilde{a}_{\text{disc}}^G(\pi)\tr \pi(\phi),
$$
where
$$
\widetilde{a}_{\text{disc}}^G(\pi)\overset{\text{df}}{=}a_{\text{disc}}^G(\pi)-\frac{1}{4}
{\sum}_{\rho:\omega_{\rho}=\omega}a_{\text{disc}}^H(\rho)\Delta_{G,H}(\rho,\pi).
$$
Similarly, for any test function $\phi_0$ on $G_0$, we get the the analogous expansion:
$$
S_{\text{disc}}^{G_0}(\phi_0)=I_{\text{disc}}^{G_0}(\phi_0)-\frac{1}{4}I_{\text{disc}}^H(\phi_0^H)={\sum}_{\Pi: \omega_{\Pi}=\omega}\widetilde{a}_{\text{disc}}^{G_0}(\Pi)\tr \Pi(\phi_0),
$$
where
$$
\widetilde{a}_{\text{disc}}^{G_0}(\Pi)\overset{\text{df}}{=}a_{\text{disc}}^{G_0}(\Pi)-\frac{1}{4}
{\sum}_{\rho:\omega_{\rho}=\omega}a_{\text{disc}}^H(\rho)\Delta_{G_0,H}(\rho,\Pi).
$$
Taken together, for any pair of matching functions $\phi \mapsto \phi^{G_0}$, we get the relation:
$$
{\sum}_{\Pi: \omega_{\Pi}=\omega}\widetilde{a}_{\text{disc}}^{G_0}(\Pi)\tr \Pi(\phi^{G_0})={\sum}_{\pi: \omega_{\pi}=\omega}\widetilde{a}_{\text{disc}}^G(\pi)\tr \pi(\phi).
$$
Now, the groups $G$ and $G_0$ are isomorphic at all finite places, so by linear independence of characters on $G(\A_F^f)$ we may isolate the finite part. That is, for every irreducible admissible representation $\tau_f$ of
$G(\A_F^f)$ we have the identity:
$$
\sum_{\Pi_{\infty}: \omega_{\Pi_{\infty}}=\omega_{\infty}}\widetilde{a}_{\text{disc}}^{G_0}(\Pi_{\infty}\otimes \tau_f)\tr \Pi_{\infty}(\phi_{\infty}^{G_0})=\sum_{\pi_{\infty}: \omega_{\pi_{\infty}}=\omega_{\infty}}\widetilde{a}_{\text{disc}}^G(\pi_{\infty}\otimes \tau_f)\tr \pi_{\infty}(\phi_{\infty}),
$$
for any pair of matching functions $\phi_{\infty} \mapsto \phi_{\infty}^{G_0}$. Of course, we assume that the fixed representation $\tau_f$ has central character $\omega_f$. In order to compare the coefficients of both sides, we invoke the character relations of Shelstad. However, first we remove the tildes by assuming $\tau_f$ is $non$-$endoscopic$, in the sense that
$$
\text{$\forall$ $\rho$ $\exists$ $v\nmid \infty$: $\Delta_{G_v,H_v}(\rho_v,\tau_v)=0$.}
$$
Furthermore, we assume $\tau_f$ is $tempered$, meaning that $\tau_v$ is tempered for almost all $v$. Then, by Proposition 6.3 on page 556 in [BR] (which assumes the congruence relation on page 555), the representation $\Pi_{\infty}\otimes \tau_f$ can only be automorphic for tempered $\Pi_{\infty}$; once we know $\Pi_{\infty}$ is cohomological for some $\xi$. It must be, if it contributes to the left-hand side of the above trace identity.

\subsection{Shelstad's character relations at infinity}

Recall our standing assumption: $\tau_f$ is a non-endoscopic tempered representation of $G(\A_F^f)$.
We wish to compare the trace identity from the previous section with the character relations of Shelstad [Sh1]. In the archimedean case, the local Langlands correspondence is completely known [La1]. Thus, the irreducible admissible representations of $G_{0,\infty}$ are partitioned into finite subsets $\Pi_{\mu}$, the $L$-packets, parameterized by $L$-parameters $\mu$. Such a $\mu$ corresponds to a tuple
$$
\text{$\mu=(\mu_v)_{v|\infty}$, $\y$ $\mu_v:W_{\R}\rightarrow {^LG_{0,v}}$.}
$$
Here $W_{\R}$ is the Weil group of $\R$, a non-trivial extension of $\{\pm 1\}$ by $\C^*$, and $\mu_v$
is a so-called admissible homomorphism. For each $v$, associated with $\mu_v$, there is an $L$-packet $\Pi_{\mu_v}$ of representations of $G_{0,v}$. Then $\Pi_{\mu}$ is the tensor product $\otimes_{v|\infty}\Pi_{\mu_v}$.
The group $G_{\infty}$ is compact mod center, so its $L$-packets are singletons, parameterized by the subset of $discrete$ parameters $\mu$; this means no $\mu_v$ maps into a parabolic subgroup. We denote the corresponding representation of $G_{\infty}$ by $\pi_{\mu}$. When $\mu_v$ is discrete, the corresponding packet for $G_{0,v}$ consists of two discrete series representations. Precisely one of them is generic (that is, it has a Whittaker model), whereas the other is holomorphic. We will use the notation:
$$
\Pi_{\mu_v}=\{\Pi_{\mu_v}^g,\Pi_{\mu_v}^h\}.
$$
Thus, the $L$-packet $\Pi_\mu$ consists of exactly $2^d$ discrete series representations.

\begin{lem}
For a fixed $\mu$, the coefficient $a_{\text{disc}}^{G_0}(\Pi_{\infty}\otimes \tau_f)$ is constant for $\Pi_{\infty}\in \Pi_{\mu}$.
\end{lem}

\noindent $Proof$. For a fixed $\tau_f$, we introduce the following distribution on $G_{0,\infty}$,
$$
T\overset{\text{df}}{=}{\sum}_{\Pi_{\infty}: \omega_{\Pi_{\infty}}=\omega_{\infty}}a_{\text{disc}}^{G_0}(\Pi_{\infty}\otimes \tau_f)\tr \Pi_{\infty}.
$$
Then $T$ is a $stable$ distribution: Suppose $\phi$ is a function on $G_{0,\infty}$ with vanishing stable orbital integrals. Then it matches the zero function on $G_{\infty}$. By the trace identity from the previous section, $T(\phi)=0$. Therefore, we can expand $T$ as
$$
\text{$T={\sum}_{\mu}c_{\mu}\tr\Pi_{\mu}$, $\y$ $\tr\Pi_{\mu}\overset{\text{df}}{=}{\sum}_{\Pi_{\infty} \in \Pi_{\mu}}\tr\Pi_{\infty}$.}
$$
Since $\tau_f$ is tempered, only the tempered $\Pi_{\infty}$ can contribute to this sum. So, we are summing over the tempered parameters $\mu$. That is, those whose image is bounded when projected onto the dual group. By linear independence of characters for $G_{0,\infty}$, we can match the coefficients and get what we aimed for,
$$
a_{\text{disc}}^{G_0}(\Pi_{\infty}\otimes \tau_f)=c_{\mu}
$$
for all $\Pi_{\infty}\in \Pi_{\mu}$. Here $\mu$ is an arbitrary (tempered) $L$-parameter. $\square$

\medskip

\noindent Now we can group together the terms on the left-hand side of the trace identity into packets, and invoke the following character relation due to Shelstad [Sh1]:
$$
\text{$\tr \Pi_{\mu}(\phi_{\infty}^{G_0})=\tr\pi_{\mu}(\phi_{\infty})$, $\y$ $\mu$ $discrete$.}
$$
Moreover, $\tr \Pi_{\mu}(\phi_{\infty}^{G_0})=0$ when $\mu$ is $not$ discrete. Again, using linear independence of characters on the group $G_{\infty}$, we deduce the equality of multiplicities:
$$
a_{\text{disc}}^G(\pi_{\mu}\otimes \tau_f)=a_{\text{disc}}^{G_0}(\Pi_{\infty}\otimes \tau_f),
$$
for all $\Pi_{\infty}\in \Pi_{\mu}$ and an arbitrary discrete $L$-parameter $\mu$. This equality immediately enables us to transfer automorphic representations from $G$ to $G_0$, and vice versa! The proof of Theorem B on Jacquet-Langlands transfer is complete.

\noindent {\sc{Department of Mathematics, Princeton University, USA.}}

\noindent {\it{E-mail address}}: {\texttt{claus@princeton.edu}}

\end{document}